\documentclass[12pt,a4paper,reqno,dvips]{article}%
\usepackage{amsfonts}
\usepackage{amsmath}
\usepackage{amssymb}
\usepackage{graphicx}
\newtheorem{theorem}{Theorem}

\newtheorem{proposition}[theorem]{Proposition}

\title{\bf Iterated Differential Forms III: Integral Calculus}
\author{\sc{A.~M.~Vinogradov}\thanks{{\bf e}-{\it mail}: \texttt{vinograd@unisa.it}} and \sc{L.~Vitagliano}\thanks{{\bf e}-{\it mail}: \texttt{luca\_vitagliano@fastwebnet.it}}\\
\small{DMI, Universit\`a degli Studi di Salerno}\\ \small{and INFN, Gruppo collegato di Salerno,}\\
\small{Via Ponte don Melillo, 84084 Fisciano (SA), Italy}}
\begin{document}
\maketitle
\begin{abstract}
 Basic elements of integral calculus over algebras of
iterated differential forms $\Lambda_k, \;k<\infty$, are
presented. In particular, defining complexes for modules of
integral forms are described and the corresponding berezinians and
complexes of integral forms are computed. Various applications and
the integral calculus over the algebra $\Lambda_{\infty}$ will be
discussed in subsequent notes.
\end{abstract}

\maketitle
\newpage
\section{Integral Calculus over Graded Algebras}

In this note we follow the notation and definitions of
\cite{vv06,vv06b} (see also \cite{v72,v81,v89,n03}). The necessary
generalities concerning integral calculus over graded commutative
algebras are collected below by following the approach of
\cite{kv98,v96}.

Let $\mathcal{G}=(G,\mu)$ be a grading group, $\Bbbk$ a field of
zero--characteristic and $A$ a $\mathcal{G}$--graded commutative algebra.
$\Lambda(A)=\bigoplus_{s}\Lambda^{s}(A)$ stands for the $\Bbbk$--algebra of
differential forms over $A$, which is naturally $\mathcal{G}\oplus%
\mathbb{Z}
$--graded. So, the de Rham complex over $A$ reads
\[
0\longrightarrow A\overset{d}{\longrightarrow}\Lambda^{1}(A)\overset
{d}{\longrightarrow}{}\cdots{}\longrightarrow\Lambda^{s}(A)\overset
{d}{\longrightarrow}{}\cdots
\]

The $A$-module of graded, skew--symmetric, $s$--derivations of $A$ with values
in an $A$--module $P$ is denoted by $\mathrm{D}_{s}(A,P)$. Put $\mathrm{D}%
_{\bullet}(A,P)=\bigoplus_{s\geq0}\mathrm{D}_{s}(A,P)$ assuming
that $\mathrm{D}_{0}(A,P)=P$. $A$--modules
$\mathrm{Diff}_{A}(P,\Lambda ^{s}(A))$, $s\geq0$, of graded linear
differential operators that send $P$ to $\Lambda^{s}(A)$ form the
complex
\begin{equation}
0\longrightarrow\mathrm{Diff}_{A}(P,A)\overset{w_{P}}{\longrightarrow
}\mathrm{Diff}_{A}(P,\Lambda^{1}(A))\overset{w_{P}}{\longrightarrow}%
\cdots\longrightarrow\mathrm{Diff}_{A}(P,\Lambda^{s}(A))\overset{w_{P}%
}{\longrightarrow}\cdots\label{wP}%
\end{equation}
with $w_{P}(\Delta)=d\circ\Delta$,
$\Delta\in\mathrm{Diff}_{A}(P,\Lambda(A))$. The cohomology
$\widehat{P}\stackrel{\mathrm{def}}{=}H(w_{P})$ of complex
(\ref{wP}) carries a natural graded $A$--module structure given by
\[
a[\Delta]=(-1)^{|a|\cdot|\Delta|}[\Delta\circ a],\quad a\in A,\,[\Delta
]\in\widehat{P},
\]
$\Delta\in\mathrm{Diff}_{A}(P,\Lambda(A))$, $w_{P}(\Delta)=0$.
$A$--module $\widehat{P}$ is called \emph{adjoint to} $P$ and
(\ref{wP}) the \emph{defining complex of} $\widehat{P}$.

A graded linear differential operator $\square:P\longrightarrow
Q$, $Q$ being an $A$--module, induces the co-chain map
\[
\widetilde{\square}:\mathrm{Diff}_{A}(Q,\Lambda(A))\ni\Delta\longmapsto
(-1)^{|\Delta|\cdot|\square|}\Delta\circ\square\in\mathrm{Diff}_{A}(P,\Lambda
(A)),
\]
i.e., $w_{Q}\circ\widetilde{\square}=\widetilde{\square}\circ w_{P}$. The
corresponding map in cohomology $\widehat{\square}:\widehat{Q}\longrightarrow
\widehat{P}$, which is a differential operator, is called the \emph{adjoint
to} $\square$ \emph{operator}.

Elements of $A$--module
$\Sigma_{s}(A)\overset{\mathrm{def}}{=}\widehat {\Lambda^{s}(A)}$
are called \emph{integral }$s$\emph{--forms} over $A$ and
$\mathcal{B}(A)\overset{\mathrm{def}}{=}\widehat{A}=\Sigma_{0}(A)$
is called the \emph{berezinian} of $A$. The module
$\Sigma(A)=\bigoplus
_{s\geq0}\Sigma_{s}(A)$ is naturally supplied with a $\mathcal{G}\oplus%
\mathbb{Z}
$--graded right $\Lambda(A)$--module structure given by
\begin{equation}
\Sigma_{l+m}(A)\times\Lambda^{l}(A)\ni(Z,\omega)\longmapsto\left\langle
Z,\omega\right\rangle \overset{\mathrm{def}}{=}[\nabla\circ\omega]\in
\Sigma_{m}(A) \label{pairing}%
\end{equation}
where $Z=[\nabla]\in\Sigma_{l+m}(A)$, $\nabla\in\mathrm{Diff}_{A}%
(\Lambda^{l+m}(A),\Lambda(A))$, $d\circ\nabla=0$.

The complex of integral forms is the adjoint to the de Rham one:
\[
0\longleftarrow\mathcal{B}(A)\overset{\widehat{d}}{\longleftarrow}\Sigma
_{1}(A)\longleftarrow\cdots{}\overset{\widehat{d}}{\longleftarrow}\Sigma
_{s}(A)\longleftarrow\cdots
\]
The following \textquotedblleft\emph{right}\textquotedblright\emph{Leibnitz
rule} holds%
\begin{equation}
\langle\widehat{d}Z,\omega\rangle=\left\langle Z,d\omega\right\rangle
+(-1)^{l}\widehat{d}\left\langle Z,\omega\right\rangle , \label{Leibnitz}%
\end{equation}
with $Z\in\Sigma(A)$, $\omega\in\Lambda^{l}(A)$.

For an $A$--module $P$ denote by
$\mathrm{Diff}_{A}^{>}(P,\Lambda(A))$ the right $A$--module
structure on the vector space of linear differential operators
acting on $P$ and with values in $\Lambda(A)$, i.e.,
\[
a^{>}\square\overset{\mathrm{def}}{=}(-1)^{|a|\cdot|\square|}\square\circ a,\quad
a\in A,\;\square\in\mathrm{Diff}_{A}^{>}(P,\Lambda(A)).
\]
There takes place a natural isomorphism (see, e.g., \cite{kv98})
\[
\mathrm{Diff}_{A}^{>}(\Lambda^{s}(A),\Lambda(A))\ni\square\longmapsto
X_{\square}\in\mathrm{D}_{s}(A,\mathrm{Diff}_{A}^{>}(A,\Lambda(A)))
\]
and a natural pairing
\begin{equation}
\mathrm{D}_{l+m}(A,\mathcal{B}(A))\times\Lambda^{l}(A)\ni(X,\omega
)\longmapsto\langle\!\langle X,\omega\rangle\!\rangle\in\mathrm{D}%
_{m}(A,\mathcal{B}(A)) \label{pairing2}%
\end{equation}
defined by
\[
\langle\!\langle X,\omega\rangle\!\rangle(a_{1},\ldots,a_{m})\overset
{\mathrm{def}}{=}i_{X}(\omega\wedge da_{1}\wedge\cdots\wedge da_{m}).
\]
$a_{1},\ldots,a_{m}\in A$. The pairing (\ref{pairing2}) supplies
$\mathrm{D}_{\bullet}(A,\mathcal{B}(A))$ with a graded right $\Lambda
(A)$--module structure.

A natural $0$--degree homomorphism of right $\Lambda(A)$--modules
\begin{equation}
\chi_{A}:\Sigma(A)\longrightarrow\mathrm{D}_{\bullet}(A,\mathcal{B}%
(A))\label{chiA}%
\end{equation}
is defined as follows. If $Z=[\nabla]\in\Sigma_{s}(A)$, $\nabla\in
\mathrm{Diff}_{A}(\Lambda^{s}(A),\Lambda(A))$, then
\[
\chi_{A}(Z)(a_{1},\ldots,a_{s})\overset{\mathrm{def}}{=}[X_{\nabla}%
(a_{1},\ldots,a_{s})]\in\mathcal{B}(A),\quad a_{1},\ldots,a_{s}\in A.
\]
Namely,
\[
\chi_{A}(\left\langle Z,\omega\right\rangle )=\langle\!\langle\chi
_{A}(X),\omega\rangle\!\rangle,\quad Z\in\Sigma(A),\;\omega\in\Lambda(A).
\]

\begin{proposition}
If $\Lambda^{1}(A)$ is a projective and finitely generated $A$--module, then
$\chi_{A}$ is an isomorphism.
\end{proposition}

In other words, in a smooth situation the isomorphism $\chi_{A}$
gives an exact description of integral forms exclusively in terms
of the berezinian $\mathcal{B}(A)$.

\section{Defining Complex of Integral Forms over Iterated Differential Form
Algebras}

In this section the complex of integral forms over the algebra of
geometric iterated differential forms on a smooth manifold is
described. This description and the results presented in the rest
of this note are generalized almost automatically to any
\textquotedblleft smooth\textquotedblright{} situation, for
instance, to super-manifolds. However, for simplicity's sake we do
not report them here. In what follows $M$ stands for a smooth
$n$--dimensional manifold, $(x^{1},\ldots,x^{n})$ for a local
chart in it and all functors of differential calculus over
commutative algebras and representing them objects are specialized
to the category of geometric modules over the algebra $C^{\infty
}(M)$ in order to be in conformity with the standard differential
geometry (see \cite{n03}). Accordingly, below
$\Lambda_{k}=\Lambda_{k}(M)$ stands for $k$ times iterated
geometric differential forms over $C^{\infty}(M)$ (see
\cite{vv06}). Put also $\Lambda_{k+1}^{s}=\Lambda^{s}(\Lambda
_{k})\subset\Lambda_{k+1}$ and note that $\Lambda_{k+1}^{s}$ is a projective
$\Lambda_{k}$--module
and $\Lambda_{k+1}^{1}$ is locally freely generated by elements $d_{k+1}%
d_{L}x^{\mu}$, \ $\mu=1,\ldots,n$, \  $L\subset\{1,\ldots,k\}$. Elements of
the dual basis are denoted by $\tfrac{\partial}{\partial d_{L}x^{\mu}}%
\in\mathrm{D}(\Lambda_{k},\Lambda_{k})\simeq\mathrm{Hom}_{\Lambda_{k}}%
(\Lambda_{k+1}^{1},\Lambda_{k})$, $\mu=1,\ldots,n$, $L\subset\{1,\ldots,k\}$,
i.e.,
\[
\tfrac{\partial}{\partial d_{L}x^{\mu}}(d_{J}x^{\nu})=\left\{
\begin{array}
[c]{ll}%
1, & \quad\text{if }\nu=\mu\text{ and }J=L\\
0, & \quad\text{otherwise}%
\end{array}
\right. .
\]

Integral geometric $p$--forms over $\Lambda_{k}$ are cohomology classes of the
complex
\[
0\longrightarrow\mathrm{Diff}_{\Lambda_{k}}(\Lambda_{k+1}^{p},\Lambda
_{k})\overset{w_{k,p}}{\longrightarrow}\mathrm{Diff}_{\Lambda_{k}}%
(\Lambda_{k+1}^{p},\Lambda_{k+1}^{1})\overset{w_{k,p}}{\longrightarrow}%
\cdots\longrightarrow\mathrm{Diff}_{\Lambda_{k}}(\Lambda_{k+1}^{p}%
,\Lambda_{k+1}^{s})\overset{w_{k,p}}{\longrightarrow}\cdots
\]
where $w_{k,p}(\Delta)=d_{k+1}\circ\Delta$, $\Delta\in\mathrm{Diff}%
_{\Lambda_{k}}(\Lambda_{k+1}^{p},\Lambda_{k+1})$. The $\Lambda_{k}$--module
$\mathrm{Diff}_{\Lambda_{k}}(\Lambda_{k},\Lambda_{k+1}^{s})$ is locally freely
generated by elements
\[
d_{k+1}d_{L_{1}}x^{\mu_{1}}\wedge\cdots\wedge d_{k+1}d_{L_{s}}x^{\mu_{s}%
}\wedge\tfrac{\partial}{\partial d_{J_{1}}x^{\nu_{1}}}\circ\cdots\circ
\tfrac{\partial}{\partial d_{J_{r}}x^{\nu_{r}}},\quad r\geq0\quad
\]
with $\mu_{1},\ldots,\mu_{s},\nu_{1},\ldots,\nu_{r}=1,\ldots,n$ and
$L_{1},\ldots,L_{s},J_{1},\ldots,J_{r}\subset\{1,\ldots,k\}$.

\section{Berezinian of the Algebra of Iterated Differential Forms and Adjoint de Rham
Differentials}

Put $\nu(k)\overset{\mathrm{def}}{=}2^{k-1}n$.
\begin{theorem}\label{main1}\quad
\begin{enumerate}
\item[(i)] $H^{s}(w_{k,0})=0$, if $s\neq\nu(k)$,
\item[(ii)] There exists a natural isomorphism of $\Lambda_{k}$--modules
$\beta_{k}:\Lambda_{k}\longrightarrow\mathcal{B}(\Lambda
_{k})=H^{\nu(k)}(w_{k,0})$.
\end{enumerate}
\end{theorem}

It is worth stressing that the definition of the Berezinian we use
in this note is different from the standard one in the theory of
super-manifolds (see, for instance, \cite{l80,b87}). Even though
these definitions are, as it can be shown, equivalent (see
\cite{v96}) the \emph{conceptuality} of the former makes it much
more preferable for our goals.

According to theorem \ref{main1} the element $\zeta_{k}=\beta_{k}%
(1_{\Lambda_{k}})$ freely generates $\Lambda_{k}$--module $\mathcal{B}%
(\Lambda_{k})$ and its local description is as follows. Let $\mathcal{U}%
=\{(x^{1},\ldots,x^{n})\}$ be a local chart and
$I_{1},\ldots,I_{2^{k-1}}$ (resp.,~$J_{1},\ldots,J_{2^{k-1}}$) be
all subsets of $\{1,\ldots,k\}$ composed of an even (resp.,~odd)
number of elements. It is assumed that each of these two families
of subsets is ordered once for all according to the subscripts.
Put
\begin{align*}
\Omega_{\mathcal{U}} &  \overset{\mathrm{def}}{=}d_{k+1}d_{I_{1}}x^{1}%
\wedge\cdots\wedge d_{k+1}d_{I_{2^{k-1}}}x^{1}\wedge\cdots\wedge
d_{k+1}d_{I_{1}}x^{n}\wedge\cdots\wedge d_{k+1}d_{I_{2^{k-1}}}x^{n}\in
\Lambda_{k+1}^{\nu(k)}(\mathcal{U}),\\
\Delta_{\mathcal{U}} &  \overset{\mathrm{def}}{=}\tfrac{\partial}{\partial
d_{J_{1}}x^{1}}\circ\cdots\circ\tfrac{\partial}{\partial d_{J_{2^{k-1}}}x^{1}%
}\circ\cdots\circ\tfrac{\partial}{\partial d_{J_{1}}x^{n}}\circ\cdots
\circ\tfrac{\partial}{\partial d_{J_{2^{k-1}}}x^{n}}\in\mathrm{Diff}%
_{\Lambda_{k}(\mathcal{U})}(\Lambda_{k}(\mathcal{U}),\Lambda_{k}%
(\mathcal{U})).
\end{align*}
Then $\zeta_{k}|_{\mathcal{U}}=[\square]\in\mathcal{B}(\Lambda_{k}%
(\mathcal{U}))$ with
\[
\square\overset{\mathrm{def}}{=}\Omega_{\mathcal{U}}\wedge\Delta_{\mathcal{U}%
}\in\mathrm{Diff}_{\Lambda_{k}(\mathcal{U})}(\Lambda_{k}(\mathcal{U}%
),\Lambda_{k+1}^{\nu(k)}(\mathcal{U})).
\]
Now define the homomorphism $\chi_{k,s}:\Sigma_{s}(\Lambda_{k}%
)\longrightarrow\mathrm{D}_{s}(\Lambda_{k},\Lambda_{k})$ by
putting
\[
(\chi_{k,s}(Z)(\omega_{1},\ldots,\omega_{s}))\,\zeta_{k}=\chi_{\Lambda_{k}%
}(Z)(\omega_{1},\ldots,\omega_{s}), \quad
Z\in\Sigma_{s}(\Lambda_{k}),
\]
(see section 1).
\begin{theorem}
$\chi_{k,s}$ is an isomorphism of $\Lambda_{k}$--modules.
\end{theorem}
So, $\Sigma_{s}(\Lambda_{k})$ is identified with
$\mathrm{D}_{s}(\Lambda _{k},\Lambda_{k})$ via $\chi_{k,s}$ and,
in view of this identification, the complex of integral forms over
$\Lambda _{k}$ reads
\[
\Lambda_{k}\overset{\widehat{d}_{k+1}}{\longleftarrow}\mathrm{D}(\Lambda
_{k},\Lambda_{k})\longleftarrow\cdots\overset{\widehat{d}_{k+1}}%
{\longleftarrow}\mathrm{D}_{s}(\Lambda_{k},\Lambda_{k})\longleftarrow\cdots .
\]
A description of the differential $\widehat{d}_{k+1}$ in these
terms is as follows.

Let $Z\in\mathrm{D}_{s}(\Lambda_{k},\Lambda_{k})$. Note that the
multi-derivation
$\widehat{d}_{k+1}Z\in\mathrm{D}_{s-1}(\Lambda_{k},\Lambda_{k})$
is completely determined by its values
$\langle\widehat{d}_{k+1}Z,\Omega\rangle$ on forms
$\Omega\in\Lambda_{k+1}^{s-1}$. But in view of (\ref{Leibnitz})
computation of these values is reduced to computation of the
action of $\widehat{d}_{k+1}$ on integral $1$--forms, i.e., on $\mathrm{D}%
(\Lambda_{k},\Lambda_{k})$. The latter is described in the
following proposition.

\begin{proposition}
If
$Z\in\Sigma_{1}(\Lambda_{k})\simeq\mathrm{D}(\Lambda_{k},\Lambda_{k})$,
then $\widehat{d}_{k+1}Z=\widehat{Z}(\zeta_{k})$ with $\widehat
{Z}:\mathcal{B}(\Lambda_{k})\simeq\Lambda_{k}\longrightarrow\mathcal{B}%
(\Lambda_{k})\simeq\Lambda_{k}$ being the adjoint to $Z$ operator.
\end{proposition}

If, locally,
\[
Z=\sum_{\mu,L}Z_{L}^{\mu}\tfrac{\partial}{\partial d_{L}x^{\mu}},\quad
Z_{L}^{\mu}\in\Lambda_{k},
\]
with the summation running over $\mu=1,\ldots,n$ and all subsets
$L$ of $\{1,\ldots ,k\}$, then
\[
\widehat{d}_{k+1}Z=-\sum_{\mu,L}(-1)^{\left\vert d_{L}\right\vert
\cdot(\left\vert d_{L}\right\vert +\left\vert Z\right\vert )}\tfrac{\partial
}{\partial d_{L}x^{\mu}}(Z_{L}^{\mu}).
\]

\begin{theorem}
$H_{i}(\widehat{d}_{k+1})\simeq H^{\nu(k)-i}(\Lambda_{1},d_{1})$.
\end{theorem}

As an example we give an explicit description of the adjoint operator $\widehat{d}%
_{l}:\mathrm{D}_{s}(\Lambda_{k},\Lambda_{k})\longrightarrow\mathrm{D}%
_{s}(\Lambda_{k},\Lambda_{k})$ to the $l$--th iterated de Rham
differential
$d_{l}:\Lambda_{k+1}^{s}\longrightarrow\Lambda_{k+1}^{s}$ for
$l\leq k$.

\begin{proposition}
Let $Z\in\mathrm{D}_{s}(\Lambda_{k},\Lambda_{k})$. Then
\begin{align*}
(\widehat{d}_{l}Z)(\omega_{1},\ldots,\omega_{s}) &  =%
{\textstyle\sum\nolimits_{l}}
(-1)^{\left\vert d_{l}\right\vert (\left\vert Z\right\vert +\left\vert
\omega_{1}\right\vert +\cdots+\left\vert \omega_{i-1}\right\vert )}%
Z(\omega_{1},\ldots,\omega_{i-1},d_{l}\omega_{i},\omega_{i+1},\ldots
,\omega_{s})\\
&  -d_{l}(Z(\omega_{1},\ldots,\omega_{s})).
\end{align*}

\end{proposition}

\section{The Trace of an Endomorphism}

Now we shall illustrate the developed theory by a simple example
showing that the notion of trace is a natural part of integral
calculus. Consider with this purpose an endomorphism of the
tangent bundle of $M$, or, equivalently, a (1,1)-tensor field on
$M$. It is naturally interpreted as a $\Lambda^{1}%
(M)$--valued derivation of the algebra $C^{\infty}(M)$ and hence
as an integral form over the algebra $\Lambda=\Lambda(M)$. Namely,
let $X\in\mathrm{D}(C^{\infty} (M),\Lambda^{1}(M))$. According to
the standard definition, $\operatorname{tr}X$, the trace  of $X$,
is the smooth function on $M$ uniquely determined by the identity
$i_{X}(\omega)=\operatorname{tr}X\;\omega$, $\omega\in
\Lambda^{n}(M)$, that can be rewritten as
\[
p_{n}\circ i_{X}=\operatorname{tr}X\;p_{n},
\]
where $p_{n}:\Lambda(M)\longrightarrow\Lambda^{n}(M)$ is a natural projection.

Now, according to the previous section $\zeta_{1}=[\square]\in\mathcal{B}%
(\Lambda)$, where $\square\in\mathrm{Diff}_{\Lambda}(\Lambda,\Lambda_{2})$ is
locally given by
\[
\square=d_{2}x^{1}\wedge\cdots\wedge d_{2}x^{n}\wedge\tfrac{\partial}{\partial
d_{1}x^{1}}\circ\cdots\circ\tfrac{\partial}{\partial d_{1}x^{n}}.
\]
The coordinate-free description of $\square$ is
\[
\square=(-1)^{n(n-1)/2}\kappa\circ p_{n},
\]
with $\kappa:\Lambda_{2}\longrightarrow\Lambda_{2}$ being the
involution that interchanges differentials $d_{1}$ and $d_{2}$
(see \cite{vv06}). Now, remembering that
$i_{X}\in\mathrm{D}(\Lambda,\Lambda)$ is an integral $1$--form
over $\Lambda$, we have:
\[
\widehat{d}_{2}(i_{X})\simeq\widehat{i_{X}}(\zeta_{1})=[\square\circ
i_{X}]=(-1)^{n(n-1)/2}[\kappa\circ p_{n}\circ i_{X}]=\operatorname{tr}%
X\;\zeta_{1}.
\]
Since $\zeta_{1}$ is a \emph{canonical} integral form, this result
allows to identify $\operatorname{tr}X$ with the integral form
$\widehat{d}_{2}(i_{X})\in\mathcal{B}(\Lambda )=\Lambda$.

\section{Final Remarks}

Our interest in integral calculus with iterated differential forms
is motivated by the problem of unification of various natural
integration procedures that one meets in differential geometry,
theoretical physics, etc. From what was said above one can see
that the developed calculus puts in common frames usual integrals,
interpreted as the de Rham cohomology, and traces of
endomorphisms. These facts give some simple arguments in favor of
iterated forms in this connection. Further results, which this
approach allows to get on in this direction, will be given in
separate publications.

\end{document}